\documentclass[11pt, a4paper,reqno]{amsart}

\usepackage{amsfonts,amsmath,amsthm}
\usepackage{amssymb}
\usepackage{latexsym}

\usepackage[english]{babel}

\theoremstyle{plain}
\newtheorem{theorem}{Theorem}[section]
\newtheorem{lemma}[theorem]{Lemma}
\newtheorem{corollary}[theorem]{Corollary}
\newtheorem{prop}[theorem]{Proposition}

\theoremstyle{definition}
\newtheorem{definition}[theorem]{Definition}

\theoremstyle{remark}
\newtheorem{remark}[theorem]{Remark}

\newcommand{\eps}{\varepsilon}
\newcommand{\N}{{\mathbb N}}

\newcommand{\R}{{\mathbb R}}
\newcommand{\ifff}{if and only if }
\newcommand{\OSM}{(\Omega,\Sigma,\mu)}

\newcommand{\Id}{\mathrm{Id}}

\newcommand{\dopu}{{:}\allowbreak\ }

\newcommand{\rest}[2]{#1\raisebox{-0.3ex}{\mbox{$\mid_{#2}$}}}

\newcommand{\loglike}[1]{\mathop{\rm #1}\nolimits}

\newcommand{\DPr}{\loglike{DPr}}
\newcommand{\ex}{\loglike{ex}}
\newcommand{\supp}{\loglike{supp}}
\newcommand{\sign}{\loglike{sign}}
\newcommand{\dist}{\loglike{dist}}
\newcommand{\lin}{\loglike{lin}}
\newcommand{\diam}{\loglike{diam}}
\newcommand{\codim}{\loglike{codim}}
\newcommand{\innt}{\loglike{int}}

\newcommand{\mytilde}[1]{\mathbin{\tilde#1}}

\newcommand{\SDt}{\loglike{{\mathcal S}\!{\mathcal D}}}
\newcommand{\narr}{\loglike{{\mathcal N}\!\!{\mathcal A}\!{\mathcal R}}}

\newcounter{abc}   % Counter f\"{u}r statements-environment wird deklariert
\newcounter{iiiii} % Counter f\"{u}r aequivalenz-environment wird deklariert

\newenvironment{aequivalenz}
{\setcounter{iiiii}{0}
\begin{list}%
{{\rm (\roman{iiiii})}}%  Falls die items nicht angegeben sind: i)u.s.w.
{\usecounter{iiiii}
\parsep=0pt plus 1pt
\topsep=1pt plus 2pt minus 1pt
\itemsep=1pt plus 2pt minus 1pt
\leftmargin=3\baselineskip \labelsep=.6\baselineskip
\labelwidth=2.4\baselineskip
\rightmargin 0pt}%
}%               Das war das zweite Argument von "newenvironment"
{\end{list}}

\newenvironment{statements}%
{\setcounter{abc}{0}
\begin{list}%
{{\rm (\alph{abc})}}%  Falls die items nicht angegeben sind: (a) u.s.w.
{\usecounter{abc}
\parsep=0pt plus 1pt
\topsep=1pt plus 2pt minus 1pt
\itemsep=1pt plus 2pt minus 1pt
\leftmargin=3\baselineskip \labelsep=.6\baselineskip
\labelwidth=2.4\baselineskip
\rightmargin 0pt}%
}%               Das war das zweite Argument von "newenvironment"
{\end{list}}

{\begin{list}%
{$\bullet$}%  Falls die items nicht angegeben sind: Punkte
{\leftmargin=3\baselineskip
\labelsep=\baselineskip \labelwidth=2.5\baselineskip
\parsep=0pt plus 1pt
\topsep=1pt plus 2pt minus 1pt
\itemsep=1pt plus 2pt minus 1pt
\rightmargin 0pt}%
}%               Das war das zweite Argument von "newenvironment"
{\end{list}}

\newcommand{\bea}{\begin{eqnarray*}}
\newcommand{\eea}{\end{eqnarray*}}
\newcommand{\beq}{\begin{equation}}
\newcommand{\eeq}{\end{equation}}
\newcommand{\begsta}{\begin{statements}}
\def\endsta{\end{statements}}
\newcommand{\begaeq}{\begin{aequivalenz}}
\def\endaeq{\end{aequivalenz}}

\newcommand{\iy}{\infty}
\def\DP{Daugavet property}
\newcommand{\pel}{Pe{\l}\-czy\'{n}\-ski}
\newcommand{\BS}{Banach space}
\newcommand{\rad}{Radon-Nikod\'ym }

\numberwithin{equation}{section}

\begin{document}

\title {Quotients of Banach spaces with the Daugavet property}

\author{Vladimir Kadets, Varvara Shepelska and Dirk Werner}

\date{May 19, 2008}

\subjclass[2000]{Primary 46B04; secondary 46B03, 46B25, 47B38}

\keywords{Daugavet property, Daugavet equation, rich subspace,
poor subspace, narrow operator, Radon-Nikod\'ym property}

 \thanks{Research of the first named author  was partially
 supported by a grant from the {\it Alexander-von-Humboldt
Foundation}. }

\address{Department of Mechanics and Mathematics,
Kharkov National University,
 pl.~Svobody~4,  61077~Kharkov, Ukraine}
\email{vova1kadets@yahoo.com; shepelskaya@yahoo.com}

\address{Department of Mathematics, Freie Universit\"at Berlin,
Arnimallee~6, \qquad {}\linebreak D-14\,195~Berlin, Germany}
\email{werner@math.fu-berlin.de}

\begin{abstract}
We consider a general concept of Daugavet property with respect to
a norming subspace. This concept covers both the usual Daugavet
property and its weak$^*$ analogue. We introduce and study analogues for
narrow operators and rich subspaces in this general setting and
apply the results to show that a quotient of $L_1[0,1]$ over an
$\ell_1$-subspace  can fail the Daugavet property.
The latter answers a question posed to us by A.~\pel\ in the negative.
\end{abstract}

\maketitle

\thispagestyle{empty}

\section{Introduction}

Throughout the paper $X$ stands for a \BS. Recall that $X$ has the
\textit{\DP}  if the identity
 \beq \label{0-eq1}
  \|\Id+T\|=1+\|T\|,
 \eeq
called {\it the Daugavet equation}, holds true for every rank-one
operator $T\dopu X \to X$.
(We shall find it convenient to abbreviate this by writing
$X\in \DPr$.)
It is known that in this case (\ref
{0-eq1}) holds for the much wider class of so-called {\it
narrow operators}. This class includes  all strong
\rad operators (which map the unit ball into a set with the \rad
property) and in particular compact and weakly compact operators,
operators not fixing a copy of $\ell_1$ and linear combinations
of the above mentioned types of operators \cite{KSW}.

Among the spaces with the \DP\ are $C(K)$-spaces and
vector-valued $C(K)$-spaces for perfect compact
Hausdorff spaces $K$, $L_1(\mu)$-
and vector-valued $L_1(\mu)$-spaces for non-atomic $\mu$, wide classes of
Banach algebras, but also some exotic spaces like Talagrand's space
from \cite{Tala5} (see \cite{KadSSW}) or Bourgain-Rosenthal's space
from \cite{BourRos} (see \cite{KW}). All the
spaces with the \DP\ are non-reflexive, moreover they cannot have
the \rad property and necessarily contain ``many" copies of
$\ell_1$ \cite{KadSSW}.

There are several results on the  stability of the \DP\ under
passing to ``big" subspaces or quotients over ``small" subspaces.
In particular, if $X\in \DPr$, then $X/E \in \DPr$ for every
reflexive subspace $E \subset X$ \cite{Shv1}. A preliminary
version of that theorem appeared in \cite{KadSSW} for $X =
L_1[0,1]$. In a private conversation after a talk by the
third-named author on the results of \cite{KadSSW}, A.~\pel\ asked
whether the last result can be generalized for $E \subset
L_1[0,1]$ being not necessarily reflexive, but having the \rad
property (RNP). The question appeared quite non-trivial for the
authors of \cite{KadSSW}, maybe because the efforts were
concentrated on attempts to prove that $L_1[0,1]/E \in \DPr$ if
$E$ has the  RNP. This question was reiterated in \cite{Shv1} and
\cite{KadKalWer}. In this paper we are now going to present a
negative answer.

Our approach to \pel's question will be indirect. If
 the answer was positive, this would mean that every subspace
$E \subset L_1[0,1]$ with the RNP has the following ``smallness"
property, since the RNP is hereditary:
$L_1[0,1]/F \in \DPr$ for every subspace $F \subset E$.
We introduce this property formally in a general setting (we call it
{\it poverty}\/) and characterise it geometrically. Then we give a
description of poor subspaces of $L_1\OSM$, and using that
description we present an $\ell_1$-subspace $E$ of $L_1[0,1]$ which
is not poor. Since $\ell_1$ has the RNP, this leads to a
counterexample.

To do all this we use duality arguments, but in order to be able to
apply these arguments we have to consider a generalisation of the \DP.
Let us recall that a subspace $Y \subset X^*$ is said to be {\it
norming}\/  (more precisely $1$-norming)
if
$$
\|x\|= \sup  \{ |y^*(x)|\dopu y^*\in Y,\ \|y^*\|\le1\}
$$
for every $x \in X$.
This is equivalent to saying that the closed unit ball of $Y$ is
weak$^*$ dense in the closed unit ball of~$X^*$.
Throughout the paper $Y$ denotes a norming
subspace of $X^*$.

\begin{definition} \label{1}
We say that $X$  has the \textit{\DP\ with respect to $Y$ $(X\in
\DPr(Y))$} if the Daugavet equation (\ref{0-eq1}) holds true for
every rank-one operator $T\dopu X \to X$ of the form $T =
y^*\otimes x$, where $x \in X$ and $y^* \in Y$.
\end{definition}

 This generalisation of the ordinary \DP\ was introduced in an
equivalent form in \cite{BKSW}. It was motivated by the fact that
the \DP\ is not stable under passing to  ultraproducts
(this was proved in \cite{KW} and an open problem at the time when
\cite{BKSW} was written),
but the ultraproduct of spaces with the \DP\ has
the \DP\ with respect to the ultraproduct of the dual spaces. The
basic motivation for us in the present paper is that the \DP\
does not generally  pass to the dual, but it is obvious that
if the original space $X$ has the \DP, then $X^* \in \DPr(X)$.

The structure of the paper is as follows. Section~2 contains
characterisations of the \DP\ with respect  to $Y$ in terms of
slices, similar to \cite{KadSSW}. Our eventual aim is to  study
``small'' subspaces $Z\subset X$ of spaces with the \DP, called
poor subspaces; this will be done in Section~5. It will turn out
that $Z\subset X$ is poor if and only if $Z^\bot \subset X^*$
enjoys a variant of the properties called richness and wealth in
\cite{KSW}. Such spaces are defined by means of a smallness
property of the corresponding quotient map called narrowness.
Narrow operators are studied in our context in Section~4. In order
to prove that weakly compact operators on dual spaces are
$X$-narrow we have included the technical Section~3 about convex
combinations of slices. Finally, in Section~6 we characterise poor
subspaces of $C(K)$ and $L_1(\mu)$ and derive that both spaces
contain copies of $\ell_1$ that are not poor provided they are
separable, $K$ is perfect and $\mu$ is atomless. This will lead to
a negative answer to \pel's question mentioned above (see
Theorem~\ref{6-theo4}).

A reader who is interested in that result only might wish to skip
Section~3 and Section~4 apart from Definition~\ref{defi3.2}, and
he or she might also wish to only consider the case $X=Y^*$ in
Section~5.

Much of the paper follows the lines of \cite{KSW}, and we omit proofs
if they don't differ much from those in that paper.

We use standard notation such as $B_{X}$  and $S_{X}$ for the unit
ball and the unit sphere of a Banach space $X$, and we employ the
notation
$$ S_U(x^{*},\eps)=
\Bigl\{x\in U\dopu x^{*}(x) > \sup_{u \in U} x^{*}(u) -\eps
\Bigr\}
$$
for the slice of a bounded convex subset $U \subset X$ determined
by $x^{*}\in S_{X^{*}}$. In the case of $U=B_{X}$ we omit the
index $U$ in the notation above:
$$
S(x^{*},\eps)= \{x\in B_{X}\dopu x^{*}(x) > 1-\eps\}.
$$
For $\eps>0$ and $x \in S_{X}$ we consider the weak$^*$ slice of the
dual ball $B_{X^{*}}$, i.e.,
$$
S(x,\eps)= \{x^{*}\in B_{X^{*}}\dopu x^{*}(x) > 1-\eps\},
$$
as a particular case of a slice. When we have a need to stress
in what space the slice is considered we use notation like $S(X,
x^{*}, \eps)$ or $S(X^{*}, x, \eps)$. The symbol $\ex C$ stands for the set
of extreme points of a set $C$. In this paper we deal with real
Banach spaces although our results extend to the complex case with
minor modifications.

\section{Basic descriptions of the generalised Daugavet property}

In this section we collect generalisations of characterisations of
the
 standard \DP\ to the setting of $\DPr(Y)$.
Let us start with a simple lemma about slices; cf.\
\cite[Lemma~2.1]{KadSSW}.

%%%%%%%%%%%%%%%%%%%%%%%%%%%%%%%
\begin{lemma}\label{1-L1}
The following statement holds true in any Banach space $X$ for any
norming subspace $Y \subset X^*$:

 Let $y\in S_{X}$, $x_0^{*}\in
S_{Y}$ and $\eps \in (0 ,1)$. Assume that there is some $x \in
S(x_0^{*},\eps/8)$ such that $\|x+y\|> 2-\eps/8$. Then there is an
$x_1^{*}\in S_{Y}$ such that $S(x_1^{*},\eps/8) \subset
S(x_0^{*},\eps)$ and $\|e+y\| > 2-\eps$ for all $e \in
S(x_1^{*},\eps/8)$.

% \begsta
%\item
%Let $y\in S_{X}$, $x_0^{*}\in S_{Y}$ and $\eps \in (0 ,1)$. Assume
%that there is some $x \in S(x_0^{*},\eps/8)$ such that $\|x+y\|>
%2-\eps/8$. Then there is an $x_1^{*}\in S_{Y}$ such that
%$S(x_1^{*},\eps/8) \subset S(x_0^{*},\eps)$ and $\|e+y\| >
%2-\eps$ for all $e \in S(x_1^{*},\eps/8)$.
%%\item
%%Let $y^{*}\in S_{X^{*}}$, $x_0\in S_{X}$ and $\eps \in (0 ,1)$.
%%Assume that there is some $x^* \in S(x_0,\eps/8)$ such that
%%$\|x^*+y^*\|> 2-\eps/8$. Then there is an $x_1\in S_{X}$ such
%%that $S(x_1,\eps/8) \subset S(x_0,\eps)$ and $\|e^*+y^*\| >
%%2-\eps$ for all $e^* \in S(x_1,\eps/8)$.
%\endsta
\end{lemma}
%%%%%%%%%%%%%%%%%%%%%%%%%%%%%%%%%%
\begin{proof}
Since $\|x+y\|> 2-\eps/8 $ there is an $x^{*}\in S_Y$ such that
$x^{*}(x+y)> 2-\eps/8$. Then
 \beq \label{1-eq3}
x^{*}(x)> 1-\frac{\eps}8 \text{\ \ and \ \ } x^{*}(y)>
1-\frac{\eps}8.
 \eeq
 Define
$x_1^{*}\in S_{X^{*}}$ by
 $$
 x_1^{*}= \frac{x_0^{*}+x^{*}}{\|x_0^{*}+x^{*}\|};
 $$
we remark that
 $$
 \|x_0^{*}+x^{*}\| \ge (x_0^{*}+x^{*})x > 2-\frac{\eps}4.
 $$
Then for every $e \in S(x_1^{*},\eps/8)$ we have
  \beq \label{1-eq02}
 (x_0^{*}+x^{*})(e) >
  \Bigl(1-\frac{\eps}8\Bigr)\Bigl(2-\frac{\eps}4\Bigr),
 \eeq
 so
 $$
 x_0^{*}(e)> 1-2\frac{\eps}8 -\frac{\eps}4 +
 \frac{\eps}8\frac{\eps}4 >1- \frac{\eps}2,
 $$
i.e., $e \in S(x_0^{*},\eps)$, and the inclusion
$S(x_1^{*},\eps/8) \subset S(x_0^{*},\eps)$ is proved. Further,
(\ref{1-eq02}) implies that $x^{*}(e)> 1- \eps/2$ which together
with (\ref{1-eq3}) means that $\|e+y\| \ge x^{*}(e + y) > 2-\eps$.
\end{proof}
%%%%%%%%%%%%%%%%%%%%%%%%%%%%%%%%%%

The following result is the analogue of \cite[Lemma~2.2]{KadSSW}.

\begin{theorem}\label{1-th2}
If $Y$ is a norming subspace of $X^*$, then the following
assertions are equivalent.
 \begaeq
 \item
$X$  has the \DP\ with respect to $Y$.
 \item
For every $x \in S_{X}$, for every $\eps > 0$, and for every
$y^*\in S_{Y}$ there is some $y \in S(y^*,\eps)$ such that
\begin{equation} \label{1-eq2}
\|x+ y\|\ge 2-\eps.
\end{equation}
 \item
For every $x \in S_{X}$, for every $\eps > 0$, and for every
$y^*\in S_{Y}$ there is a slice $S(y_1^*,\eps_1) \subset
S(y^*,\eps)$ with $y_1^*\in S_{Y}$ such that {\rm(\ref{1-eq2})} holds
for every $y \in S(y^*,\eps_1)$.
% \item
%For every $x^*\in S_{Y}$, for every $\eps > 0$, and for every
%slice $S(x,\eps)$ of the dual ball $B_{X^{*}}$ there is some
%$y^*\in S(x,\eps)$ such that $\|x^{*}+ y^*\|\ge 2-\eps$.
% \item
%For every $x^*\in S_{Y}$, for every $\eps > 0$, and for every
%slice $S(x,\eps)$ of the dual ball $B_{X^{*}}$ there is another
%slice $S(x_1,\eps_1) \subset S(x,\eps)$ such that $\|x^{*}+
%y^*\|\ge 2-\eps$ for every $y^*\in S(x_1,\eps_1)$ .
 \endaeq
\end{theorem}

%%%%%%%%%%%%%%%%%%%%%%%%%%%%%%%%%%
\begin{proof}
The implication (iii) $\Rightarrow$ (ii)
%and (v) $\Rightarrow$ (iv)
is obvious; the implication (ii) $\Rightarrow$ (iii)
%and (iv) $\Rightarrow$ (v)
follows from Lemma~\ref{1-L1}. What remains to prove is the
equivalence (i) $\Leftrightarrow$ (ii).
%and (i) $\Leftrightarrow$ (iv).

Let us start with (i) $\Rightarrow$ (ii). Fix some $x\in
S_{X}$, $x^* \in S_Y$ and $\eps > 0$ and consider the operator
$T\dopu X \to X$, $T := x^*\otimes x$. According to (\ref{0-eq1}),
 $\|\Id+T\|=2$, so there is a  $y\in
S_{X}$ such that $\|y + Ty\|\ge  2-\eps/2$ and $x^*(y) \ge 0$.
Substituting the value of $Ty$ we obtain that
 $$
 \|y + x^*(y)x\|\ge  2-\eps/2
 $$
which means that $x^*(y) \ge 1 - \eps/2$ (i.e., $y\in S(x^*,\eps)$)
and
 $$
 \|y + x\| \ge \|y + x^*(y)x\| - |x^*(y) - 1| \ge  2-\eps,
 $$
which proves the implication.

For the converse implication (ii) $\Rightarrow$ (i) consider an
operator $T\dopu X \to X$, $T = x^*\otimes x$, where $x \in X$ and
$x^* \in Y$. Since the validity of (\ref{0-eq1}) for $T$ implies
(\ref{0-eq1}) for all operators of the form $aT$ with $a > 0$, it
is sufficient to consider the case of $\|T\|=1$, and the
representation $T = x^*\otimes x$ can be taken in such a way that
$x \in S_X$ and $x^* \in S_Y$. Due to (ii), for every $\eps > 0$
there is a $y\in S(x^*,\eps)$ satisfying  (\ref{1-eq2}). Then
\bea
\|\Id+T\| &\ge& \|y+Ty\| = \|y + x^*(y)x\| \\
&\ge& \|y + x\| -  |x^*(y) - 1| \ge (2-\eps) - \eps,
\eea
which by arbitrariness of $\eps$ means that $\|\Id+T\|=2$.

%The equivalence (i) $\Leftrightarrow$ (iv) can be proved in the same
%way as the equivalence (i) $\Leftrightarrow$ (ii) if one rewrites
%(\ref{0-eq1}) in the equivalent form $ \|\Id_{X^*}+T^*\| = 1 +
%\|T^*\|$.
\end{proof}

%%%%%%%%%%%%%%%%%%%%%%%%%%%%%%%%%%
%%%%%%%%%%%%%%%%%%%%%%%%%%%%%%%%%%
\section{A useful tool: convex combinations of slices}

This section deals with a technical device that will be useful
in the proof of Theorem~\ref{3-theo 1}.

\begin{definition} \label{ccs1}
Let $Y$ be a norming subspace of $X^*$, and let $U \subset X$ be
convex and bounded. A subset $V \subset U$ is called a
quasi-$\sigma_U(X,Y)$ neighbourhood if it is a finite convex
combination of slices of  $U$  generated by elements of $Y$;
i.e., there are $\lambda _k \ge 0$, $ k=1,\dotsc ,n$, with $\sum_{k=1}^n
\lambda _k = 1$ and slices $S_1, \dotsc ,S_{n} \subset U$
generated by elements of $Y$ such that $ \lambda _1S_1 + \cdots +
\lambda _nS_n = V$.
\end{definition}

The following lemma is known for the ordinary weak topology (see
\cite[Lemme~5.3]{Bour-Paris}; it was rediscovered in \cite{Shv1}).
The $\sigma(X,Y)$-version proof coincides almost word-to-word with
the original one.

\begin{lemma}\label{ccs2}
Under the conditions of the above definition every relatively
$\sigma(X,Y)$-open subset $A \subset U$ contains a
quasi-$\sigma_U(X,Y)$ neighbourhood.
\end{lemma}

The next theorem and its corollary  were essentially proved by
Shvidkoy \cite{Shv1}. He considered the ordinary \DP, but the
proof in the general case is virtually the same.

\begin{theorem}\label{ccs3}
Let $Y$ be a norming subspace of $X^*$ and $X \in \DPr(Y)$. Then for
every $\varepsilon > 0$, every $x \in S_{X}$ and every
quasi-$\sigma_{B_X}(X,Y)$ neighbourhood $V$ there exists an element
$v \in V$ such that $\|v+x\|\ge 2 - \varepsilon$.
\end{theorem}

\begin{proof}
Let $V = \sum_{k=1}^n \lambda_k S_k $ be a representation of $V$
as a convex combination of slices. Using repeatedly (ii) of
Theorem~\ref{1-th2} one can construct $x_k \in S_k$ such that
$\|(x + \sum_{j < k}\lambda_j x_j) +  \lambda_k x_k\| \ge \|x +
\sum_{j < k}\lambda_j x_j\| +  \lambda_k - \eps/n$. Then $v =
\sum_{k=1}^n \lambda_k x_k$ will be the element of $V$ we need.
\end{proof}

\begin{corollary}\label{1-th3}
If $Y$ is a norming subspace of $X^*$ and $X$  has the \DP\ with
respect to $Y$, then the following is true:
For every $x \in S_{X}$, for every $\eps > 0$, and for every
$\sigma(X,Y)$-open subset $U \subset X$ intersecting $B_X$ there
is some $y \in U \cap B_X$ such that $\|x+ y\|\ge 2-\eps$.
\end{corollary}

%%%%%%%%%%%%%%%%%%%%%%%%%%%%%%%%%%
%%%%%%%%%%%%%%%%%%%%%%%%%%%%%%%%%%

\section{Narrow  operators with respect to a norming subspace}

We will eventually study subspaces satisfying a certain smallness
condition called ``poverty''; this will be dual to the notion of a ``rich''
subspace from \cite{KSW}. The latter class is defined by the
requirement that the canonical quotient map is ``narrow''. This
section deals with such operators.

First we will
recall and modify some definitions from \cite{KSW}. Let $X$, $E$ be
Banach spaces.

\begin{definition} \label{defi3.1}
An  operator $T \in L(X,E)$ is said to be a {\em strong Daugavet
operator}\/ if for every two elements $x, y \in S_{X}$ and for
every $\varepsilon > 0$ there is an element $z \in S_{X}$
such that $\|z+x\|> 2 - \varepsilon$ and $\|Tz-Ty\|< \eps$. We
denote the class of all strong Daugavet operators on $X$ by
$\SDt(X)$.
\end{definition}

Corollary~\ref{1-th3} shows that if $X \in \DPr(Y)$ then every $T
\in L(X,E)$ of the form $T = f\otimes e$, where $e \in E$ and $f
\in Y$, is a strong Daugavet operator, and conversely, thanks to
Theorem~\ref{1-th2}, if every $f \in Y \subset X^*=L(X,\R)$  is
strongly Daugavet, then $X$ has the Daugavet property with respect
to $Y$.

There is an obvious connection between strong Daugavet operators
and the Daugavet equation (cf.\ \cite[Lemma 3.2]{KSW}).

\begin{lemma} \label{3-lem1}
If $T\dopu X\to X$ is a strong Daugavet operator, then $T$
satisfies the Daugavet equation {\rm(\ref{0-eq1})}.
\end{lemma}

\begin{definition} \label{defi3.2}
Let $X \in \DPr(Y)$. An operator $T \in L(X,E)$ is said to be {
\it narrow with respect to $Y$}\/ (or \textit{$Y$-narrow} for short)
if  for every $x, e \in
S_{X}$, $\varepsilon > 0$  and every slice $S \subset B_X$
generated by an element of $Y$ and containing $e$ there is an
element $v \in S$ such that $\|x+v\|>2-\varepsilon$ and
$\|Tv-Te\|< \varepsilon$. We denote the class of all $Y$-narrow
operators on $X$ by  $\narr_Y(X)$.
\end{definition}

The notations  $\SDt(X)$ and $\narr_Y(X)$ do not mention the range
space $E$ because  the corresponding definitions do not actually
depend on the values of $T$, but only on the norms of those values,
i.e., these are not properties of the operator $T$
itself, but just of the seminorm $x\mapsto \|T(x)\|$ on $X$. For more about
this ideology see \cite{KSW}.

The following statement is a complete analogue of
\cite[Lemma~3.10(a)]{KSW},  so we omit the proof.

\begin{lemma}\label{Lccs}
Let $T\in \narr_Y(X)$. Let $S_1, \dotsc ,S_{n} \subset B_X$ be a
finite collection of slices generated by elements of $Y$, and let
$U \subset B_{X}$ be a convex combination of these slices, i.e.,
there are $\lambda _k \ge 0$ with $ \sum_{k=1}^n
\lambda _k = 1$ such that $ \lambda _1S_1 + \cdots + \lambda
_nS_n = U$. Then for every $\varepsilon > 0$, every $x_1 \in
S_{X}$ and  every $w\in U$ there exists an element $u \in U$ such
that $\|u+x_1\|> 2 - \varepsilon$ and $\|T(w-u)\|<\varepsilon$.
\end{lemma}

Let us recall an operation with operators that was introduced in
\cite{KSW}. For operators $T_{1}\dopu X \to E_{1} $ and
$T_{2}\dopu X \to E_{2} $ define
$$
T_{1}\mytilde+ T_{2}\dopu X \to E_{1}\oplus_{1}E_{2}, \quad x
\mapsto (T_{1}x, T_{2}x);
$$
i.e.,
$$
\|(T_{1} \mytilde+ T_{2})x\| = \|T_{1}x\| + \|T_{2}x\|.
$$

\begin{remark} \label{rem-defi3.2}
Let $X$, $E$ be Banach spaces, $Y \subset X^*$ be a norming
subspace, $T \in L(X,E)$. If $T \mytilde+ y^* \in \SDt(X)$ for
every $y^* \in Y$, then $T \in \narr_Y(X)$.
In the setting of $Y=X^*$ this was actually
given as the definition of a narrow
operator in \cite{KSW}, and our definition was given as an equivalent
condition in Lemma~3.10 of~\cite{KSW}.
\end{remark}

We are now going to introduce a class of operators that turn out to be
$Y$-narrow; they correspond to the strong \rad operators in
the case $Y=X^*$, which contain the weakly compact operators.
We need two technical definitions.

\begin{definition} \label{defi3.3}
Let $X$, $E$ be Banach spaces, $F \subset E^*$ be a norming subspace and
$\eps > 0$. A point $e$ of a convex subset $A \subset E$ is said to
be an {\em $(F, \eps)$-denting point}\/ if there is a functional $f \in
S_F$ and a $\delta > 0$ such that $\|e - a\| < \eps$ whenever $a
\in A$  satisfies the condition $f(a) > f(e) - \delta$.
 We say that $A \subset E$ is {\em $F$-dentable}\/ if for every $\eps > 0$
the set $A$ is contained in the closed convex hull of its $(F,
\eps)$-denting points. An operator $T \in L(X,E)$ is said to be
{\em $F$-dentable}\/ if $T(B_X)$ is $F$-dentable. An $E^*$-dentable
operator is called {\em dentable}.
\end{definition}

\begin{definition}\label{3-def 4}
An operator $T \in L(X,E)$ is said to be {\em hereditarily $F$-dentable}\/
if for every  $x^* \in X^*$ the operator
 $$
T \mytilde+ x^* \dopu X \to E \oplus_{1} \R, \quad x \mapsto (Tx,
x^*x)
 $$
 is $\tilde F$-dentable, where $\tilde F$ consists of all
 functionals $(f, \beta) \dopu E \oplus_{1} \R \to \R$, with $f \in F$ and
 $\beta \in \R $, of the form  $(f, \beta)\left((e, t)\right)= f(e) + \beta
 t$.
\end{definition}

 Remark that every strong \rad operator and in particular
every weakly compact operator is hereditarily dentable, by well-known
geometric characterisations of sets with the RNP~\cite[Chap.~3]{Bou}.

\begin{theorem}\label{3-theo 1}
If $X \in \DPr(Y)$, $T \dopu X \to E$ is a hereditarily $F$-dentable
operator and $T^*(F) \subset Y$, then $T$ is $Y$-narrow.
\end{theorem}

\begin{proof}
According to Remark~\ref{rem-defi3.2} it is sufficient to prove
that $\tilde{T}=T\tilde{+}x^{*} \in \SDt(X)$ for every $x^{*} \in
Y$. Fix $x, z \in S_{X}$ and $\eps
> 0$. By Definition~\ref{defi3.1}, to prove the theorem we have
to find an element $v \in B_X$ such that $\|x+v\|>2-\varepsilon$
and $\|\tilde{T}(z-v)\|< \varepsilon$.

Since $\tilde{T}(B_X)$ is $\tilde F$-dentable there are $\lambda
_k \ge 0$ with $\sum_{k=1}^n \lambda_k = 1$ and
$(\tilde F, \eps/2)$-denting points  $e_1, \ldots, e_n \in T(B_X)$
such that
$$
\biggl\|Tz - \sum_{k=1}^n \lambda_k e_k \biggr\| < \eps/2.
$$
By the
definition of an $\tilde F$-denting point there are slices
$S_k=S_{\tilde T(B_X)}(f_k, \eps_k)$   of $T(B_X)$ generated by elements
of $\tilde F$ such that $e_k \in S_k$ and the diameter of each of the
$S_k$ is less than $\eps/2$. Denote $W:= \sum_{k=1}^n \lambda_k
S_k$. Since $\dist(Tz, W) < \eps/2$ and $\diam W < \eps/2$ we have
\begin{equation} \label{eq3-th1-1}
\|Tz - w\| < \eps \qquad \text{for every } w \in W.
\end{equation}
Denote $y_k^*:= \tilde{T}^*f_k$. By assumption,  $y_k^* \in Y$.
Consider slices $V_k = \{v \in B_X \dopu y_k^*(x) > \eps_k\}$ and
the quasi-$\sigma_{B_X}(X,Y)$ neighbourhood $V := \sum_{k=1}^n
\lambda_k V_k$. Since $\tilde{T}(V_k) \subset S_k$ and
consequently $\tilde{T}(V) \subset W$, (\ref{eq3-th1-1}) implies
that
$$
\|Tz - Tv\| < \eps \qquad \text{for every } v \in V.
$$
It remains to apply Theorem~\ref{ccs3} to get a $v \in V$ with
$\|x+v\|>2-\varepsilon$.
\end{proof}

\begin{corollary} \label{cor5.9}
A  weak$^*$-weakly continuous
operator on a dual space $X^*\in \DPr(X)$ is
$X$-narrow.
\end{corollary}

\begin{proof}
It is clear that such an operator is weakly compact.
The hereditary dentability of a weakly compact
operator $T\dopu X^*\to E$ has
been mentioned in the above remark; it remains to observe that
$T^*(E^*)\subset X$ if $T$ is weak$^*$-weakly continuous.
\end{proof}

Note that a weakly compact adjoint operator $S^*\dopu X^*\to V^*$ is
weak$^*$-weakly continuous, i.e., $S^{**}(V^{**}) \subset X$; cf.\
\cite[Section~VI.4]{DS1}.

%%%%%%%%%%%%%%%%%%%%%%%%%%%%%%%%%%

\section{Rich and poor subspaces}

In \cite{KSW} we introduced rich subspaces $Z\subset X$, building on
work in \cite{KadPop} and \cite{PliPop}. We showed that this condition
is equivalent to saying that every superspace $Z\subset \tilde
Z\subset X$ has the \DP; the latter property was called wealth in
\cite{KSW}. We now extend and dualise these ideas.

\subsection{Richness}

The next proposition  shows a kind of stability of the
\DP\  when one passes from the original space to a ``big"
subspace.

%%%%%%%%%%%%%%%%%%%%%%%55%%%%%%%%

\begin{lemma}\label{4-L1}
Let $X\in \DPr(Y)$.  Then for every $x\in S_{X}$, every
$\varepsilon>0$ and for every separable subspace $V \subset Y$
there is an $x^*\in S_{X^*}$ such that $x^*(x)\ge1-\varepsilon$
and $\|x^*+f\|=1+\|f\|$ for all $f\in V$.
\end{lemma}

%%%%%%%%%%%%%%%%%%%%%%%55%%%%%%%%
\begin{proof}
Consider a dense  sequence
$(f_n)_{n=1}^{\infty} \subset V$ such  that every element is
repeated  infinitely many times in the sequence. Applying (v) of
Theorem~\ref{1-th2} to the slice $S(x,\eps)$ of $B_{X^*}$ and to
$f_1$ and then applying it step-by-step to $f_n$ and to the slices
obtained in the previous steps, we construct a sequence of closed slices
$\overline S(x,\eps) \supset \overline S(x_1,\eps_1) \supset
\overline S(x_2,\eps_2) \supset
\ldots$  with $\eps_n < 1/n$  such that
$\|x^*+f_n\|\ge 2-\varepsilon_{n-1}$ for all $x^*\in
S(x_n,\varepsilon_n)$. By $w^*$-compactness of all
$\overline S(x_n,\varepsilon_n)$, there is a point
$x^*\in\bigcap_{n=1}^{\infty} \overline S(x_n,\varepsilon_n)\subset
\overline S(x,\varepsilon)$. This is exactly the point we need.
\end{proof}

\begin{prop}\label{4-pr1}
Let $X\in \DPr(Y)$, and let $Z \subset X$ be a subspace such that
$Z^{\bot}$ is a separable subspace of $Y$. Then $Z\in \DPr(\rest{Y}{Z})$.
\end{prop}

\begin{proof}
Let $z \in S(Z)$,
and let
 $$
 S=\{z^* \in Z^* = X^*/Z^{\bot}\dopu \|z^*\|\le 1,\  z^*(z) \ge 1-\varepsilon\}
 $$
be a slice of $B_{Z^*}$. Fix a $[z^*]\in S(Y/Z^{\bot})$. We have
to prove the existence of  an $[x^*]\in S$ such that
$\|x^*+z^*\|=2$. Applying Lemma~\ref{4-L1} with $x=z$ and
$V=\lin(\{z^*\}\cup Z^{\bot})$ we obtain an $x^*\in S_{X^*}$
such that $x^*(z)\ge1-\varepsilon$ and
 $$
 \|x^*+f\|=1+\|f\|\qquad {\textrm{for  all }}  f\in V.
 $$
Then $[x^*]\in S$ and
 \bea
 \|[x^*+z^*]\| &=&
  \inf_{f\in Z^{\bot}} \|x^*+z^*+f\| \\
  &=&
 \inf_{f\in Z^{\bot}} (1+\|z^*+f\|) =
  1+\|[z^*]\|=2.
 \eea
\end{proof}

If a subspace $Z \subset X$ satisfies the conditions of the proposition
above then so do all the subspaces of $X$ containing $Z$.
Hence $Z$ has the property that $\tilde Z \in
\DPr(\rest{Y}{{\tilde Z}})$ for every subspace $\tilde Z \subset X$
containing $Z$. Let us formalise this property.

\begin{definition}\label{12}
Let $X\in \DPr(Y)$. A subspace $Z \subset X$ is said to be {\it
wealthy with respect to}\/ $Y$ if  $\tilde Z \in \DPr(\rest{Y}{{\tilde
  Z}})$
for every subspace $\tilde Z \subset X$ containing $Z$.
\end{definition}

Thus Proposition~\ref{4-pr1} can be rephrased by saying that $Z\subset
X$ is wealthy with respect to $Y$ if $Z^\bot$ is separable and $X\in
\DPr(Y)$.

The main result of this subsection is a characterisation of
$Y$-wealthy subspaces through $Y$-narrow operators, analogous to
\cite[Theorem 5.12]{KSW}.

\begin{definition}\label{13}
Let $X\in \DPr(Y)$. A subspace $Z\subset X$ is said to be {\em rich
with respect to}\/ $Y$ if  the quotient map $q\dopu  X
\rightarrow X/Z$ is a $Y$-narrow operator.
\end{definition}

It turns out that the following theorem holds.

\begin{theorem}\label{4-teo1}
Let $X$ be a Banach space and $Y$ be a norming subspace of $X^*$
such that $X\in \DPr(Y)$. Then for a subspace $Z\subset X$ the
following properties are equivalent:
\begaeq
\item $Z$ is wealthy with respect to $Y$.
\item $Z$ is rich with respect to $Y$.
\endaeq
\end{theorem}

The proof is very similar to the proof of \cite[Theorem~5.12]{KSW}.

\subsection{Poverty as a dual property to richness}

\begin{definition}\label{5.2-def2}
Let $X \in \DPr$. A subspace $Z \subset X$ is said to be {\it
poor}\/ if  $ X / \tilde Z \in \DPr$ for every subspace $\tilde Z
\subset Z$.
\end{definition}

Our study of poor subspaces uses duality, so let us start with a
very simple observation that we state as a proposition for easy reference.

\begin{prop}\label{5.2-prop1}
A \BS\ $X$ has the \DP\  \ifff $X^* \in \DPr(X)$. Hence, a
subspace $Z $ of a space $X$ with the \DP\  is poor \ifff
for every subspace $\tilde{Z}\subset Z$ its dual
$(X/\tilde{Z})^*=\tilde{Z}^{\bot}$ has the \DP\  with respect to
$X/\tilde Z$.
\end{prop}

Now we are ready to give the basic characterisations of poverty.

\begin{theorem}\label{5.2-theor1}
Let $X \in \DPr$. For a subspace $Z \subset X$ the following
 conditions are equivalent.
 \begaeq
\item
$Z$ is poor.
\item
$ X / \tilde Z \in \DPr$ for every subspace $\tilde{Z}\subset Z$
of codimension $\codim_Z \tilde{Z} \le 2$.
\item
$Z^{\bot}$ is a subspace of $X^*$ that is rich with respect to~$X$.
\item
For every $x^*, e^*  \in S_{X^*}$,
$\varepsilon > 0$ and for every $x \in S_{X}$ such that $e^*(x) >
1 - \eps$ there is an element $v^* \in B_{X^*}$ with the following
properties: $v^*(x) > 1 - \eps$, $\|x^*+v^*\|>2-\varepsilon$ and
$\|\rest{(e^*-v^*)}{Z}\|< \eps$; that is,
 the quotient map from $X^*$ onto $X^*/Z^{\bot}$ is narrow
with respect to $X$.
 \endaeq
\end{theorem}

\begin{proof}
(i) $\Rightarrow$ (ii) follows immediately from the definition of
poor subspaces.

Let us prove (ii) $\Rightarrow$ (i). According to
Proposition~\ref{5.2-prop1}, we have to prove that for every subspace
$Z_1\subset Z$, its dual $Z_1^{\bot}$ has the \DP\ with respect to
$X/Z_1$. Fix $Z_1\subset Z$. Applying (ii) of Theorem~\ref{1-th2} we
see that for every $x^*\in S_{Z_1^{\bot}}$, $\eps>0$ and every
$[x]\in S_{X/Z_1}$ we have to find $y^*\in S_{Z_1^{\bot}}$ such
that $y^*([x])\ge 1-\eps$ and $\|x^*+y^*\|\ge 2-\eps$.
Since $[x]\in S_{X/Z_1}$, there exists $z^*\in
S_{Z_1^{\bot}}$ such that $z^*([x])=1$. Denote
$\tilde{Z}=Z\cap \ker x^*\cap \ker z^*$. Evidently, $\tilde{Z}$ is a
subspace of $Z$ of $\codim_Z \tilde{Z}\le 2$ and $Z_1\subset
\tilde{Z}$. Also remark that
$$
1=\|[x]_{X/Z_1}\|\ge \|[x]\|_{X/\tilde{Z}}\ge z^*([x])=1,
$$
which implies $[x]\in
S_{X/\tilde{Z}}$. By our assumption $\tilde{Z}^{\bot}$ has the
Daugavet property with respect to $X/\tilde Z$, and hence for $x^*\in
S_{\tilde{Z}^{\bot}}$ and $[x]\in S_{X/\tilde{Z}}$ there is
$y^*\in S_{\tilde{Z}^{\bot}}$ such that $y^*([x])\ge 1-\eps$ and
$\|x^*+y^*\|\ge 2-\eps$. Then $y^*\in S_{\tilde{Z}^{\bot}}\subset
S_{Z_1^{\bot}}$, and it meets all the requirements.

Now we will prove that (ii) $\Leftrightarrow$ (iii).
Theorem~\ref{4-teo1} implies that
(iii) holds if and only if $Z^{\bot}$ is a subspace of $X^*$ that is
  wealthy with respect to $X$; and this is equivalent to the claim that for
every $x^*,y^*\in S_{X^*}$ the space $W
=\lin (Z^{\bot} \cup\{x^*,y^*\})$ has the Daugavet property with
respect to $X/W_\bot$ (cf.\ \cite[Lemma 5.6(iii)]{KSW}. But for a
space $\hat{Z}\supset Z^{\bot}$ the existence of $x^*,y^*\in
S_{X^*}$ such that $W=\lin (Z^{\bot} \cup\{x^*,y^*\})$ is
equivalent to the existence of a space $\tilde Z\subset Z$ such
that $W=\tilde{Z}^{\bot}$ and $\codim_Z \tilde{Z}\le 2$.
Thus we get that (iii) is equivalent to the claim that
$\tilde{Z}^{\bot}\in\DPr(X/\tilde Z)$ for every subspace $\tilde{Z}\subset
Z$ of $\codim_Z \tilde{Z}\le 2$, which is equivalent to (ii)
according to Proposition~\ref{5.2-prop1}.

The remaining equivalence (iii) $\Leftrightarrow$ (iv) is just a
reformulation of the definition of a rich subspace.
\end{proof}

As a corollary we can give a proof of the following theorem of Shvidkoy
\cite{Shv1}.

\begin{corollary}\label{5.2-cor1}
Let $X \in \DPr$ and let $Z$ be a reflexive subspace of $X$. Then
the quotient space $X/Z$ also has the \DP.
\end{corollary}

\begin{proof}
Since every subspace of a reflexive space is also reflexive, the
statement of this corollary is equivalent to the claim that every
reflexive subspace $Z$ of $X \in \DPr$ is poor. According to
Theorem~\ref{5.2-theor1}, it is sufficient to prove that
$Z^{\bot}$ is a subspace of $X^*$ that is rich with respect to $X$, i.e.,
that the quotient map $q\dopu X^*\rightarrow X^*/Z^{\bot}$ is an
$X$-narrow operator. As $X^*/Z^{\bot}$ is isometric to $Z^*$,
which is reflexive, this follows from Corollary~\ref{cor5.9}.
\end{proof}

\section{Applications to the geometry of $C(K)$ and $L_1$ }

For a compact Hausdorff space
$K$ denote by $M(K)$ the dual space of $C(K)$, i.e.,
$M(K)$ is the \BS\ of all (not necessarily positive) finite
regular Borel signed measures on~$K$. (In the sequel, all measures on
$K$ will be tacitly assumed to be finite regular Borel measures.)
We are going to prove a theorem which gives a
characterisation of  operators on $M(K)$ that are narrow
with respect to
$C(K)$. For this theorem  we will need  the following lemma in which
  $\partial A$ denotes the boundary of a set
$A \subset K$.

\begin{lemma}\label{6-lem}
Let $K$ be compact, $f\in C(K)$, and $\mu$ be some positive measure on
$K$.  Then for every $\eps>0$ there exists a step function
$\tilde{f}=\sum_{k=1}^n \beta_k \chi_{A_k}$ on $K$ such that
 $\mu(\partial A_k)=0$ for $k=1,\dots,n$, $A_1 \cup \dots \cup A_n=K$ and
$\|f-\tilde{f}\|_{\infty}<\eps$.
\end{lemma}

\begin{proof}
Since the image measure $\nu=\mu\circ f^{-1}$ on $\R$ has at most
countably many atoms, it is possible to cover $f(K)$ by finitely
many half-open intervals $I_k= (\beta_{k-1}, \beta_k]$ of length
${<\eps}$ such that $\nu(\{\beta_0,\dots,\beta_n\})=0$. Let $A_k=
f^{-1}(I_k)$; then $\tilde{f}=\sum_{k=1}^n \beta_k \chi_{A_k}$
works.
\end{proof}

\begin{theorem}\label{6-theo1}
Let $K$ be a perfect compact Hausdorff space.
An operator $T$ on  $M(K)$ is narrow
with respect to $C(K)$ \ifff for every open subset  $U \subset K$,
for every two  probability  measures $\pi_1,\pi_2$
on $U$  and for every $\eps > 0$ there is a  probability
measure $\nu$ on $U$ such that $\|T(\nu-\pi_1)\|< \eps $ and
$\|\pi_2-\nu\|>2-\eps$.
\end{theorem}

\begin{proof}
We first prove the ``only if" part. By the definition of a narrow
operator (Definition~\ref{defi3.2}), for every  $x, e  \in
S_{M(K)}$, $\varepsilon > 0$ and every weak$^*$ slice $S$ of
$B_{M(K)}$ containing  $e$, there exists  $v
\in S$ such that $\|x+v\|>2-\varepsilon$ and $\|T(e-v)\|<
\varepsilon$. Fix $\varepsilon_1>0$ and let $x=-\pi_2$ and
$e=\pi_1$. Since $U$ is open and $ \pi_1(U)=1$, we can
find $f\in C(K)$  taking values in $[0,1]$ with
$\supp f\subset U$ and $\int f\,
d\pi_1>1-\varepsilon_1$. By Definition~\ref{defi3.2}
there exists $\tilde{\nu}\in
S_{M(K)}$ such that the following inequalities hold:
$$
\int f \,d\tilde{\nu}>1-\varepsilon_1, \quad
\|T(\tilde{\nu}-\pi_1)\|<\varepsilon_1, \quad
 \|\pi_2-\tilde{\nu}\|>2-\varepsilon_1.
$$

Let $\hat{\nu}=\rest{\tilde{\nu}^{+}}{U}$. Using the properties of $f$
we have $\|\tilde{\nu}-\hat{\nu}\|<2\varepsilon_1$ and thus
$$
1-3\varepsilon_1<\|\hat{\nu}\|\le1+2\varepsilon_1, \quad
\|T(\hat{\nu}-\pi_1)\|<\varepsilon_1(1+2\|T\|), \quad
\|\pi_2-\hat{\nu}\|>2-3\varepsilon_1.
$$
Hence for $\nu={\hat{\nu}}/{\|\hat{\nu}\|}$ we have
$\|\nu-\hat{\nu}\|=\bigl| 1-\|\hat{\nu}\| \bigr|<3\varepsilon_1$ and
consequently
$$
\|\pi_2-\nu\|\ge\|\pi_2-\hat{\nu}\|-\|\nu-\hat{\nu}\|>
2-3\varepsilon_1-3\varepsilon_1=2-6\varepsilon_1,
$$
and
$$\|T(\nu-\pi_1)\|\le \|T(\hat{\nu}-\pi_1)\|+\|T(\hat{\nu}-\nu)\|<
(1+5\|T\|)\varepsilon_1.
$$
Then taking
$\varepsilon_1=\min\{\frac{\varepsilon}{6},\frac{\varepsilon}{1+5\|T\|}\}$
completes the proof of the ``only if" part.

Now consider the ``if" part. Given $\mu_1,\mu_2  \in S_{M(K)}$,
$\varepsilon > 0$ and a weak$^*$ slice $S$ of $B_{M(K)}$
containing $\mu_1$, we have to find $\nu \in S$ such that
$\|\mu_2+\nu\|>2-\varepsilon$ and $\|T(\mu_1-\nu)\|< \varepsilon$.
Since one can wiggle the slice $S$ a bit, there is, by
Lemma~\ref{6-lem}, no loss of generality in replacing $S$ by a
slice generated by a function of the form $f=\sum_{k=1}^n \beta_k
\chi_{A_k}$, where $A_1,\ldots,A_n$ are measurable sets with
$(|\mu_1|+|\mu_2|)(\bigcup_{k=1}^n \partial A_k)=0$. (Note that in
general this new slice will not be relatively weak$^*$ open.) On
the other hand, using the Hahn decomposition theorem we have
$K=\bigcup_{i=1}^{4}B_i$, where $B_1$ is a set on which $\mu_1$ is
positive and $\mu_2$ is negative, $B_2$ is a set on which $\mu_2$
is positive and $\mu_1$ is negative, and $B_3$ (resp.\ $B_4$) is a
set where  both $\mu_1$ and $\mu_2$ are positive (resp.\
negative).

Fix  $\varepsilon_1>0$ and let $G_1$ be an open set such that
$G_1\supset B_1$ and $|\mu_i|(G_1\setminus B_1)<\varepsilon_1$
($i=1,2$). Define $C_k=G_1\cap A_k$ and let $U_k= \innt C_k$,
$k=1,\ldots,n$. Clearly $C_k\setminus U_k\subset \partial A_k$, so the
$U_k$ are open  sets with the following properties: $U_k\subset
C_k$ and $(|\mu_1|+|\mu_2|)(C_k\setminus U_k)=0$.

Consider those $U_k$ for which $\mu_1(U_k\cap B_1)\ne 0$,
$\mu_2(U_k\cap B_1)\ne 0$ and define two probability measures on
$U_k$ by
$$
\mu_{i,k}=\frac{\rest{\mu_i}{U_k\cap B_1}}{\mu_i(U_k\cap B_1)}
\quad (i=1,2).
$$
By assumption there exists a probability  measure
$\hat{\nu}_k$ on $U_k$ such that
$$
\|T(\hat{\nu}_k-\mu_{1,k})\|<\eps_1
\ \text{ and }\
\|\mu_{2,k}-\hat{\nu}_k\|>2-\eps_1.
$$
Define
$\nu_k=\mu_1(U_k\cap B_1)\cdot \hat{\nu}_k$. Then we have
\begin{equation}\label{6-eq1}
\|\nu_k\|=\nu_k(U_k)=\mu_1(U_k\cap B_1),\quad
\|\rest{\mu_1}{U_k}\|-\eps_1\le\|\nu_k\|\le\|\rest{\mu_1}{U_k}\|+\eps_1
\end{equation}
and
\begin{eqnarray}\label{6-eq2}
\|\rest{\mu_2}{U_k}+\nu_k\| &=&\|\mu_2(U_k\cap B_1)\cdot
\mu_{2,k}+ \rest{\mu_2}{U_k\setminus B_1}
    +\mu_1(U_k\cap B_1)\cdot \hat{\nu}_k\| \nonumber\\
    &\ge &
\||\mu_2(U_k\cap B_1)|\cdot \mu_{2,k}-|\mu_1(U_k\cap B_1)|\cdot
\hat{\nu}_k\|- \|\rest{\mu_2}{U_k\setminus
B_1}\|\nonumber\\
     &\ge & |\mu_2|(U_k)+|\mu_1|(U_k)-4\eps_1
\end{eqnarray}
and
\begin{eqnarray}\label{6-eq3}
\|T(\nu_k-\rest{\mu_1}{U_k})\| &\le & \|T(\mu_1(U_k\cap
B_1)\cdot(\hat{\nu}_k-\mu_{1,k}))\|+ \|T(\rest{\mu_1}{U_k\setminus
B_1})\| \nonumber\\
    &\le& \eps_1(1+\|T\|).
\end{eqnarray}
For  $U_k$ with $\mu_1(U_k\cap B_1)=0$ or $\mu_2(U_k\cap B_1)=0$,
the inequalities (\ref{6-eq1})--(\ref{6-eq3}) hold with
 $\nu_k=\rest{\mu_1}{U_k\cap B_1}$.

Now define the measure $\mu_1^1$ by
$$
\rest{\mu_1^1}{U_k}=\nu_k,\quad \rest{\mu_1^1}{K\setminus
\bigcup_{k=1}^{n}U_k}=\rest{\mu_1}{K\setminus \bigcup_{k=1}^{n}U_k}.
$$
 From
(\ref{6-eq1}), (\ref{6-eq2}), and (\ref{6-eq3}) we obtain the
following properties of $\mu_1^1$:
\begin{equation}\label{6-eq4}
\|\mu_1\|-n\eps_1\le\|\mu_1^1\|\le\|\mu_1\|+n\eps_1,\quad
\biggl|\int f \,d\mu_1^1
    -\int f \,d\mu_1\biggr| \le n\varepsilon_1
\end{equation}
and
\begin{eqnarray}\label{6-eq5}
\|\rest{\mu_2}{G_1}+\rest{\mu_1^1}{G_1}\| &\ge&
\biggl\|\sum_{k=1}^n(\rest{\mu_2}{U_k}+\nu_k) \biggr\|
    -(|\mu_2|+|\mu_1|)
         \biggl(\bigcup_{k=1}^{n}C_k\setminus
         U_k \biggr) \nonumber\\
    &\ge& \sum_{k=1}^n (|\mu_2|(U_k)+|\mu_1|(U_k))-4n\eps_1\nonumber\\
    &\ge& |\mu_1|(G_1)+|\mu_2|(G_1)-(4n+2)\eps_1   , \\
%\end{eqnarray}
%and
%\begin{equation}\label{6-eq6}
\|T(\mu_1^1-\mu_1)\|
&=&
\biggl\|\sum_{k=1}^n
T(\nu_k-\rest{\mu_1}{U_k})\biggr\|\le(n+n\|T\|)\eps_1. \label{6-eq6}
%\end{equation}
\end{eqnarray}

Now define $\tilde{B}_2=B_2\setminus G_1$. Notice that
$\tilde{B}_2$ is a set of negativity for $\mu_1^1$ and a set
of positivity for $\mu_2$. Following the same lines as above we
define $G_2\supset \tilde{B}_2$ and construct $\mu_1^2\in M(K)$
such that
$$
|\mu_2|(G_2\setminus \tilde{B}_2)<\eps_1,\quad |\mu_1^1|(G_2\setminus
\tilde{B}_2)<\eps_1,\quad \|\mu_1^1\|-n\eps_1\le\|\mu_1^2\|\le
\|\mu_1^1\|+n\eps_1
$$
and
\bea
\biggl|\int f \,d\mu_1^2-\int f \,d\mu_1^1\biggr|
&\le&
 n\varepsilon_1, \\
\|T(\mu_1^2-\mu_1^1)\|
&\le&
(n+n\|T\|)\eps_1, \\
\|\rest{\mu_2}{G_2}+\rest{\mu_1^2}{G_2}\|
&\ge&
|\mu_1|(G_2)+|\mu_2|(G_2)-(4n+2)\eps_1.
\eea

 From (\ref{6-eq4}), (\ref{6-eq5}), (\ref{6-eq6}) and the above
inequalities we obtain the estimates
$$
1-2n\eps_1\le\|\mu_1^2\|\le1-2n\eps_1,\quad \biggl|\int f \,
d\mu_1^2-\int f \,d\mu_1\biggr|
    \le 2n\varepsilon_1
$$
and
\begin{eqnarray}
\|T(\mu_1^2-\mu_1)\| &\le& (2n+2n\|T\|)\eps_1,  \label{t}  \\
\|\rest{\mu_2}{G_1\cup G_2}+\rest{\mu_1^2}{G_1\cup G_2}\|
    &\ge&
|\mu_1|(G_1\cup G_2)+|\mu_2|(G_1\cup G_2)-(8n+10)\eps_1. \nonumber
\end{eqnarray}

Finally, the definition of the sets $B_3$ and $B_4$ implies that
\begin{equation} \label{mu2}
\|\mu_2+\mu_1^2\|\ge \|\mu_1\|+\|\mu_2\|-(8n+10)\eps_1=2-(8n+10)\eps_1.
\end{equation}
Hence for  $\eps_1$ small enough, the normalized signed measure
$\nu={\mu_1^2}/{\|\mu_1^2\|}$ satisfies all the required
conditions, which completes the proof of the
theorem.
\end{proof}

Applying this theorem to the operator $\mu\mapsto \rest{\mu}{Z}$
yields by Theorem~\ref{5.2-theor1}:

\begin{corollary}\label{6-cor2-}
Let $K$ be a perfect compact. A subspace $Z \subset C(K)$ is poor
\ifff  for every open subset  $U \subset K$, for every two
probability  measures $\pi_1,\pi_2$ on $U$  and for every
$\eps > 0$ there is a  probability measure $\nu$ on
$U$ such that $\|\nu - \pi_1\|_{Z^*}< \eps$ and
$\|\pi_2-\nu\|>2-\eps$.
\end{corollary}

For a closed subset $K_1$ of $K$ denote by $R_{K_1}$ the natural
restriction operator $R_{K_1}\dopu C(K) \to C(K_1)$. Note that for an
operator $S\dopu E\to F$ between Banach spaces the following
assertions are equivalent, by the (proof of) the open mapping theorem:
(i)~$S$ is onto; (ii)~$S(B_E)$ is not nowhere dense; (iii)~$0$~is an
interior point of $S(B_E)$.

\begin{corollary}\label{6-cor2}
Let $K$ be a perfect compact Hausdorff space,
$K_1 \subset K$ be a closed subset
with non-empty interior, and let $Z$ be a poor subspace of $C(K)$.
Then $R_{K_1}(B_Z)$ is nowhere dense in $B_{C(K_1)}$.
\end{corollary}

\begin{proof}
Apply Corollary~\ref{6-cor2-} with $U=\innt K_1$, $\pi_1=\pi _2$
and a sufficiently small $\eps>0$ to see that $R_{K_1}(B_Z)$
cannot contain a ball $rB_{C(K_1)}$ of radius ${r>0}$.
\end{proof}

We now deal with poor subspaces of $L_1$.
Let $\OSM$ be a finite measure space. Denote by $\Sigma^+$ the
collection of all $A \in \Sigma$ with $\mu(A) > 0$.

\begin{theorem}\label{6-theo2}
Let $\OSM$ be a non-atomic finite measure space.
An operator $T$
on $L_\iy:=L_\iy\OSM$ is narrow with respect to $L_1:=L_1\OSM$
\ifff for every $\Delta \in \Sigma^+$ and for every $\eps > 0$
there is $g \in S_{L_\iy}$ such that $g=0$ off $\Delta$
and $\|Tg\|< \eps $.
Moreover,  in the statement above $g$ can be selected non-negative.
\end{theorem}

 \begin{proof}
First we prove the ``if" part. By the definition of a narrow
operator (Definition~\ref{defi3.2})  for every $x,y\in
S_{L_\iy}$, every $f\in S_{L_1}$ such that $\int f\cdot y\,
d\lambda > 1-\delta$ (i.e., $y\in S(f,\delta)$) and every $\eps>0$ we
have to find $z\in S(f,\delta)$ such that $\|x+z\|>2-\eps$ and
$\|T(y-z)\|<\eps$.
By density of step functions we may assume
without loss of generality that
there is a partition $A_1,\dots,A_n$ of $\Omega$ such that the
restrictions of $x,y$ and $f$ to $A_k$ are constants, say $a_k$,
$b_k$ and $c_k$ respectively. Fix some $\eps_1>0$. Since
$\|x\|=1$, there exists $k$ such that
$|a_k|>1-\eps_1$. Let $B\in \Sigma^{+}$ be a subset of $A_k$ with
$\mu(B)\le \eps_1$ and $A_k\setminus B\in\Sigma^{+}$. By our
assumption there exists $\hat{z}\in S_{L_\iy}$ such that $z\ge0$,
$z$ is supported on $ B$ and $\|T(z)\|\le \eps_1$. Denote
$\tilde{z}=y+(\sign(a_1)-b_1)\hat{z}$. It is easy to see that
$\|\tilde{z}\|=1$, $\|x+\tilde{z}\|=2-\varepsilon_1$,
$\|T(y-\tilde{z})\|<\varepsilon_1$ and $\tilde{z}\in
S(f,\delta-\varepsilon_1)$. To finish the proof of this part
it is sufficient to repeat the reasoning from the end of the proof
of Theorem~\ref{6-theo1}.

Now we consider the ``only if" part. Since $T$ is narrow with
respect to $L_1$, $T$ is also a strong Daugavet operator.
Hence  as in \cite[Theorem~3.5]{KSW} we can get a
function $\tilde{g}\in S_{L_\iy}$ which satisfies all the
requirements, except being non-negative. To fix this we argue as
in \cite[Lemma~1.4]{KadPop} and finally get some non-negative
$g$ possessing all the properties listed above.
\end{proof}

Again, specialising to the restriction operator $g\in L_\infty =
(L_1)^* \mapsto \rest gZ \in Z^*$ we obtain the following
characterisation of poor subspaces.

\begin{corollary}\label{6-cor4-}
Let $\OSM$ be a non-atomic finite measure space. A subspace $Z
\subset  L_1\OSM$ is poor \ifff for every $\Delta \in \Sigma^+$
and for every $\eps > 0$ there is $g \in S_{L_\iy}$ such that
$ g=0$ off $\Delta$ and $\|g\|_{Z^*}< \eps $. Moreover,
in the statement above $g$ can be selected non-negative.
\end{corollary}

For a subset $A \in \Sigma^+$ denote by $Q_{A}$ the natural
restriction operator $Q_{A}\dopu L_1\OSM \to L_1(A, \rest{\Sigma}{A},
\rest{\mu}{A})$.

\begin{corollary}\label{6-cor4}
Let $\OSM$ be a non-atomic finite measure space, $A \in \Sigma^+$
and let $Z$ be a poor subspace of $L_1\OSM$. Then $Q_{A}(B_Z)$  is
nowhere dense in $B_{L_1(A, \Sigma|_A, \mu)}$.
\end{corollary}

\begin{proof}
Apply Corollary~\ref{6-cor4-}  with $\Delta=A$
and a
sufficiently small $\eps>0$ to see that $Q_A(B_Z)$ cannot contain
a ball $rB_{L_1(A)}$ of radius ${r>0}$.
\end{proof}

The Corollaries~\ref{6-cor2} and~\ref{6-cor4} look very similar. The
next definition extracts the  significant  common feature.

\begin{definition}\label{6-def2}
Let $X \in \DPr$. A subspace $E \subset X$ is said to be a {\it
bank}\/ if $E$ contains an isomorphic copy of $\ell_1$ and for every
poor subspace $Z$ of $X$,  $q_E(B_Z)$ is nowhere dense in $B_{X/E}$
(here $q_E$ denotes the natural quotient map $q_E \dopu X \to X/E$).
If $E \subset X$ is a bank, then $B_{X/E}$ will be called  the
{\it asset}\/ of $E$.
\end{definition}

In this terminology a poor subspace cannot cover a ``significant
part" of a bank's asset.

\begin{theorem}\label{6-theo3}
Let $X \in \DPr$ and $E \subset X$ be a bank with separable asset.
Then $X$ contains a copy of $\ell_1$ which is not poor in $X$.
\end{theorem}

\begin{proof}
Let $\{e_n\}_{n \in \N} \subset \frac12 B_E$ be equivalent to the
canonical basis of $\ell_1$ and let $\{x_n\}_{n \in \N} \subset B_E$
be  a sequence such that $\{q_E(x_n)\}_{n \in \N}$ is dense in
$B_{X/E}$. Then, if one selects a sufficiently small $\eps > 0$,
the sequence of $u_n = e_n + \eps x_n \in B_E$ is
still equivalent to the canonical basis of $\ell_1$, and the image of
this sequence under  $q_E$ equals $\{\eps q_E(x_n)\}_{n \in \N}$, which
is dense in $\eps B_{X/E}$. This means that the closed linear span
of $\{u_n\}_{n \in \N}$ is the copy of $\ell_1$ we need.
\end{proof}

The next theorem is an immediate corollary of
Theorem~\ref{6-theo3}.

\begin{theorem}\label{6-theo4}
In every $C(K)$-space with perfect metric compact $K$ and in every
separable $L_1\OSM$-space with non-atomic $\mu$ there is a
subspace  isomorphic to $\ell_1$ that  is not poor.
\end{theorem}

\begin{proof}
Corollary~\ref{6-cor2} implies that if $K$ is a perfect compact and $K_1
\subset K$ is a proper closed subset with non-empty interior, then
$C_0(K \setminus K_1):=\{f \in C(K) \dopu f(t)=0$ $ \forall t \in
K_1\}$ is a bank with $B_{C(K_1)}$ being its asset.
Corollary~\ref{6-cor4} implies that if $\OSM$ is a non-atomic finite measure
space and $A \in \Sigma^+$, then $L_1(\Omega \setminus A)$ is a bank
with $B_{L_1(A)}$ being its asset. Separability of these assets
follows from the separability of the spaces $C(K)$ and $L_1\OSM$
considered.
It is left to apply Theorem~\ref{6-theo3}.
\end{proof}

Theorem~\ref{6-theo4} answers \pel's question mentioned in the
introduction in the negative since it provides a non-poor
$\ell_1$-subspace  $Z\subset L_1[0,1]$. By definition this means that
for some subspace $\tilde Z\subset Z$,  $L_1[0,1]/\tilde Z$
fails the \DP; but $Z$ has
the RNP and so does its subspace~$\tilde Z$. In fact, by
Theorem~\ref{5.2-theor1} one can choose $\tilde Z$ of
codimension~${\le2}$, hence $\tilde Z$ is isomorphic to $\ell_1$ as
well. Let us some up these considerations.

\begin{corollary}\label{6-cor5}
There is a subspace $E\subset L_1[0,1]$ that is isomorphic to $\ell_1$
and hence has the RNP, but $L_1[0,1]/E $  fails the \DP.
\end{corollary}

%%%%%%%%%%%%%%%%%%%%%%%%%%%%%%%%%%%%%%%%%%%%%%%%%%%%%%%%%%%%%%%%%%%

\section{Some open questiions}

1. Is it true that every separable space with the \DP\ has an
$\ell_1$-subspace which is not poor?

2. Can the separability condition in
Theorem~\ref{6-theo3}  be omitted?

3. Is it true that every subspace without copies of $\ell_1$ of a
space with the \DP\ is poor? We don't even know the answer  in the
case of $C[0,1]$.

4. Is it true that if $X \in \DPr$ and $Y \subset X$ is a subspace
with a separable dual, then the quotient space $X/Y$ also has the
\DP? This question also appears in \cite{Shv1}.

%%%%%%%%%%%%%%%%%%%%%%%%%%%%%%%%%%%%%%%%%%%%%%%%%%%%%%%%%%%%%%%%%%%%%
%%%%%%%%%%%%%%%%%%%%%%%%%%%%%%%%%%%%%%%%%%%%%%%%%%%%%%%%%%%%%%%%%%%%%

\end{document}

%%%%%%%%%%%%%%%%%%%%%%%%%%%%%%%%%%%%%%%%%%%%%%%%%%%%%%%%%%%%%%%%%%%%%
%%%%%%%%%%%%%%%%%%%%%%%%%%%%%%%%%%%%%%%%%%%%%%%%%%%%%%%%%%%%%%%%%%%%%